\documentclass[12pt]{article}
\usepackage{amsfonts,amsmath}
\usepackage{amssymb,latexsym}
\usepackage[dvips]{epsfig}
\usepackage{amscd}
\usepackage{psfrag}
\usepackage[all,cmtip]{xy} 

\title{How to categorify one-half of quantum $\mathfrak{gl}(1|2)$}
\author{ Mikhail Khovanov}
\date{July 25, 2010}

\newtheorem{lemma}{Lemma}
\newcommand{\oplusop}[1]{{\mathop{\oplus}\limits_{#1}}}
 
 \newcommand{\sumoop}[2]{{\mathop{\sum}\limits_{#1}^{#2}}}

\begin{document} 

\maketitle
\baselineskip 14pt
 
\def\R{\mathbb R}
\def\Q{\mathbb Q}
\def\Z{\mathbb Z}
\def\N{\mathbb N} 
\def\C{\mathbb C}
\def\l{\lbrace}
\def\r{\rbrace}
\def\o{\otimes}
\def\lra{\longrightarrow}
\def\Hom{\mathrm{Hom}}
\def\RHom{\mathrm{RHom}}
\def\Id{\mathrm{Id}}
\def\mc{\mathcal}
\def\mf{\mathfrak} 
\def\Ext{\mathrm{Ext}}
\newcommand{\define}{\stackrel{\mbox{\scriptsize{def}}}{=}}
\newcommand{\ii}{{\bf i}}
\newcommand{\jj}{{\bf j}}
\def\drawing#1{\begin{center}\epsfig{file=#1}\end{center}}
 \def\yesnocases#1#2#3#4{\left\{
\begin{array}{ll} #1 & #2 \\ #3 & #4
\end{array} \right. }
\newcommand{\LOT}{H^-} 

\begin{abstract} 
We describe a collection of differential graded rings that categorify weight 
spaces of the positive half of the quantized universal enveloping algebra 
of the Lie superalgebra $\mf{gl}(1|2)$. 
\end{abstract}


\psfrag{X1}{$a, b\in \{1,2\}$} 
\psfrag{X2}{$\delta_{a,2}\delta_{b,1}\delta_{c,2}$} 
\psfrag{X3}{$e_{m_1}$}\psfrag{X4}{$e_{m_2}$}\psfrag{X5}{$e_{m_k}$}
\psfrag{X6}{$1_{1 2^{(m)}}= $} \psfrag{X7}{$1_{2^{(m)}1}=$} \psfrag{X8}{$e_m$} 
\psfrag{X9}{$1_{2^k 1 2^{m-k}}=$}
\psfrag{X10}{$1_{1 2^{(m_1)}1\dots 1 2^{(m_k)} 1 }=$} 
\psfrag{X11}{$\tau(\alpha_a')$} \psfrag{X12}{$\tau(\beta_a')$}  
\psfrag{X13}{$\tau(\alpha_a'')$} \psfrag{X14}{$\tau(\beta_a'')$} 
\psfrag{X15}{}
 \psfrag{X16}{{\Large $y_k = \ \ \  \sumoop{a=1}{r'(1)}$}}
\psfrag{X17}{{\Large $\sumoop{a=1}{r''(1)}$}}
\psfrag{X18}{{\Large $1_{2^k 1 2^{m-k}}= \ \ \sumoop{a=1}{r'(1)}$}} 
\psfrag{X19}{{\Large $X = $}}
\psfrag{X20}{{\Large $dX = $}}
\psfrag{X21}{{$j_1$}}
\psfrag{X22}{{$j_2$}}
\psfrag{X23}{{$k_1$}}
\psfrag{X24}{{$k_2$}}
\psfrag{X25}{{$(j_1,j_2)$-crossing}}
\psfrag{X26}{{$(k_1,k_2)$-crossing}}
\psfrag{X27}{{$e_m$}}


\noindent 
{\bf Lie superalgebra $\mf{gl}(1|2)$, the positive half and its quantum version} 


\noindent 
The Lie superalgebra $\mf{gl}(n|m)$ is defined by partitioning $(n+m)\times (n+m)$ 
matrices into 4 blocks: diagonal blocks of size $n\times n$ and $m\times m$ and 
off-diagonal blocks of size $n\times m$ and $m\times n$. Matrices with nonzero entries 
only in the diagonal, respectively off-diagonal blocks, are called \emph{even}, 
respectively \emph{odd}. Elementary matrix $E_{ij}$ 
is even if $i,j\le n$ or $i,j\ge n+1$ and odd otherwise.

The superbracket $[A,B]$ of matrices is defined as the 
usual bracket $AB-BA$ if at least one of $A$ and $B$ is even and as 
anticommutator $AB+BA$  if both $A$ and $B$ are odd. With these conventions, 
the superbracket satisfies the super analogues~\cite{Kac1}  of the antisymmetry 
and the Jacobi identity: 
\begin{eqnarray} 
   [a,b] & = & - (-1)^{p(a)p(b) } [b,a] ,    \\                
 {} [a,[b,c]] & = & [[a,b],c] + (-1)^{p(a)p(b)} [b,[a,c]],  
\end{eqnarray} 
where $p(a)=0$ (resp. $1$) if $a$ is even (resp. odd). 

The universal enveloping algebra $UL$ of a Lie superalgebra $L$ is defined in the 
same way as for Lie algebras. If $L=L_0\oplus L_1$ is the decomposition of $L$ 
into the sum of its even and odd parts, $UL$ can be identified, as a vector space, 
with $S(L_0)\otimes \Lambda(L_1)$, the tensor product of the symmetric algebra of 
$L_0$ and exterior algebra of $L_1$, once bases of $L_0$ and $L_1$ are fixed. 
It is a Hopf algebra in the category of super vector 
spaces, with $\Delta(x)=x\otimes 1 + 1\otimes x$ for $x\in L$. 

The decomposition of $\mf{gl}(n)$ into 
the direct sum of strictly upper-triangular, diagonal, and strictly lower-triangular matrices 
generalizes to $\mf{gl}(n|m)$. We denote by $\mf{gl}(n|m)^+$ the Lie superalgebra 
of strictly upper-triangular $(n|m)$-matrices and by $U^+(n|m)$ its universal 
enveloping algebra (say, over $\Q$). 

Note that $U^+(1|1)\cong \Lambda(E_{12})$, the exterior algebra on 
one generator $E_{12}$. In particular, $U^+(1|1)$ is two-dimensional, 
with basis $\{1,E_{12}\}$, and $E_{12}^2=0$.   

The universal enveloping algebra $U^+(1|2)$ (we also denote it $U^+$) 
has two generators $E_1:=E_{12}$ and $E_2:=E_{23}$, 
first odd, second even, and defining relations 
\begin{eqnarray*} 
   E_1^2 & = & 0 , \\
   2 E_2 E_1 E_2 & = & E_1 E_2^2 + E_2^2 E_1. 
\end{eqnarray*} 
The Lie superalgebra $\mf{gl}(1|2)^+$ has a basis $\{E_2, E_1, [E_1,E_2]\}$. 
The first of these generators is even, the other two are odd, so that $U^+$ 
has a basis $\{E_2^m E_1^{\epsilon}[E_1,E_2]^{\epsilon'}\}$ where $m\in \Z_+$ 
and $\epsilon, \epsilon'\in \{0,1\}$. The set $\{E_1^{\epsilon}E_2^m E_1^{\epsilon'}\}$ 
is also a basis, with $m,\epsilon,\epsilon'$ in the same range as above, except that 
$m\not= 0$ if $\epsilon=\epsilon'=1$. 
 
The quantum deformation $U_q^+=U_q^+(1|2)$ of $U^+$ is a $\Q(q)$-algebra with 
the same generators as $U^+$ and defining relations 
\begin{eqnarray} 
   E_1^2 & = & 0 , \label{eq-one} \\
   {} [2] E_2 E_1 E_2 & = & E_1 E_2^2 + E_2^2 E_1. \label{eq-two} 
\end{eqnarray} The deformation simply transforms coefficient $2$ in the second relation to quantum 
$[2]=q+q^{-1}$. The latter relation can be rewritten 
\begin{equation} 
    E_2 E_1 E_2 = E_1 E_2^{(2)} + E_2^{(2)} E_1, \quad \quad 
E_2^{(2)}:=\frac{E_2^2}{[2]}, \label{eq-three}
\end{equation} 
where $E^{(m)}:= \frac{E^m}{[m]!}$ denotes $m$-th quantum divided 
power of $E$. 

Equipped with comultiplication 
$$ \Delta(E_i) = E_i \otimes 1 + 1 \otimes E_i, \quad \quad i=1,2,$$ 
$U^+_q$ becomes a twisted bialgebra, in the sense of~\cite[Chapter 1]{Lus1}, in the 
category of super vector spaces.  

Define the integral form $U^+_{\Z}=U^+_{\Z}(1|2)$ 
to be the $\Z[q,q^{-1}]$-subalgebra of $U^+_q$ generated by $E_1$ and divided powers 
$E_2^{(m)}$ over all $m\ge 0$. As a free $\Z[q,q^{-1}]$-module, 
$U^+_{\Z}$ has a basis $\{E_1^{\epsilon} E_2^{(m)}E_1^{\epsilon'}\}$, 
with $m\in \Z_+,$ $\epsilon, \epsilon'\in \{0,1\}$, and $m\not= 0$ if $\epsilon=\epsilon'=1$.  The set of defining relations in $U^+_{\Z}$ can be taken to be 
\begin{eqnarray} 
 E_1^2 & = & 0    \\
 E_2^{(k)} E_2^{(m-k)} & = & 
\left[ \hspace{-0.05in} \begin{array}{c} m \\ k \end{array} \hspace{-0.05in}\right]  
E_2^{(m)},  
\quad \quad 0 \le k \le m,  \label{eq-four} \\
E_2^{(k)} E_1 E_2^{(m-k)} & = & 
\left[ \hspace{-0.05in} \begin{array}{c} m-1 \\ k \end{array}  \hspace{-0.05in} \right]  
E_1 E_2^{(m)}  + 
\left[ \hspace{-0.05in} \begin{array}{c} m-1 \\ k-1 \end{array} \hspace{-0.05in}  \right] 
E_2^{(m)} E_1 ,  \   0 < k < m, \label{eq-five} 
\end{eqnarray} 
where square brackets denote quantum binomials. 
Notice that $E_2^{(k)} E_1 E_2^{(m-k)}$ is a linear combination 
of $E_1 E_2^{(m)}$ and $E_2^{(m)}E_1$ with coefficients in $\Z_+[q,q^{-1}]$. 
The weight space of $U^+_q$ containing these products is naturally isomorphic 
to the corresponding weight space of $U^+_q(\mf{sl}(3))$, since the only relation 
that contributes to the size of this weight space is (\ref{eq-two}) in both algebras. 
Moreover,  $\{  E_1 E_2^{(m)}, E_2^{(m)}E_1\}$ is the canonical 
basis of this weight space of quantum $\mf{sl}(3)$, see \cite[Example 14.5.4]{Lus1}.  

\vspace{0.15in} 

 
\noindent 
{\bf Categorification of positive half of quantum $\mf{sl}(2)$ } 


\noindent 
We recall categorification of $U_q^+(\mf{sl}(2))$ following \cite{KL1, Lau1}. 
Fix a ground field $\Bbbk$.  
The nilHecke algebra $H_m$ is the algebra of endomorphisms of 
$\Bbbk[x_1, \dots, x_m]$ generated by multiplication by $x_i$ endomorphisms 
(also denoted $x_i$) and divided difference operators 
$$ \partial_i(f) = \frac{f-s_i(f)}{x_i -x_{i+1}} ,$$ 
where $s_i(f)$ is $f$ with $x_i, x_{i+1}$ transposed. We depict the identity endomorphism 
of $\Bbbk[x_1, \dots, x_m]$ by $m$ vertical lines, multiplication by $x_i$ endomorphism 
by the dot on the $i$-th strand counting from left, and $\partial_i$ by the $i$-th crossing: 

\drawing{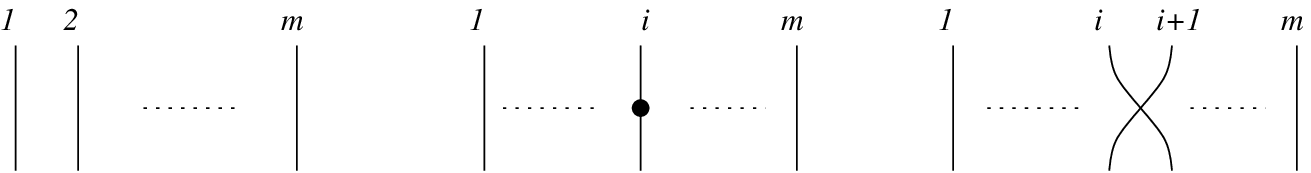} 

The following is a  set of defining relations for $H_m$:  
\begin{eqnarray*} 
& & x_i x_j = x_j x_i,  \\
& &  x_i \partial_j = \partial_j x_i, \ \ \ i\not= j,j+1, \quad \quad   
  \partial_i \partial_j = \partial_j \partial_i, \ \ \  |i - j|>1, \\
& & \partial_i^2 = 0, \quad \quad \partial_i \partial_{i+1}\partial_i = \partial_i \partial_{i+1} 
\partial_i , \\ 
& & x_i \partial_i - \partial_i x_{i+1} = 1, \quad 
      \partial_i x_i - x_{i+1}\partial_i = 1. 
\end{eqnarray*} 
The defining relations say that far away dots and crossings can be isotoped past  
each other and, in addition, the following diagrammatic equalities hold 

\drawing{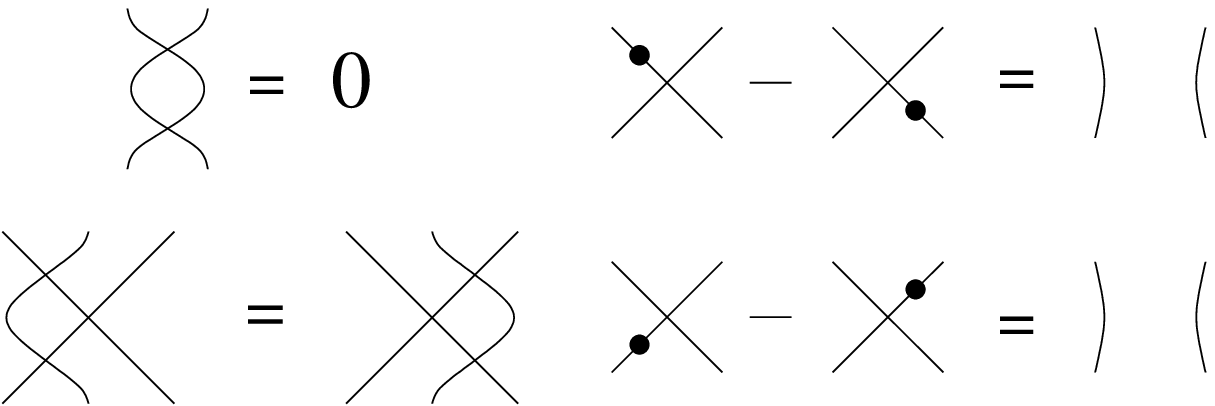} 

The center of $H_m$ is isomorphic to the ring of symmetric polynomials 
in $x_1, \dots, x_m$, and $H_m$ is isomorphic to the algebra of $m!\times m!$ 
matrices with coefficients in the center $Z(H_m)$. The minimal idempotent 
$e_m\in H_m$, given for $m=3$ by the diagram below, 

\drawing{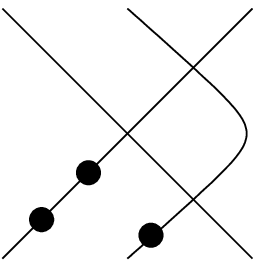} 

\noindent 
provides a Morita equivalence between $H_m$ and its center, via the 
bimodules $H_m e_m$ and $e_m H_m$, the first an $(H_m, Z(H_m))$-bimodule, 
the second a $(Z(H_m), H_m)$-bimodule. Idempotent $e_m$ is the product  
of the maximal permutation word in divided difference operators and 
$x_1^{m-1} x_2^{m-2}\dots x_{m-1}$. 

Algebra $H_m$ is graded, with $\deg(x_i)=2$ and $\deg(\partial_i)=-2$, 
so that $\deg(e_m)=0$, and the above Morita equivalence is that of graded rings.  

The Grothendieck group $K_0(A)$ of a $\Z$-graded associative ring $A$ is a
$\Z[q,q^{-1}]$-module with generators $[P]$, 
over finitely-generated graded projective $A$-modules $P$, and defining relations 
$[P]=[P']+[P'']$ whenever $P\cong P'\oplus P''$
and $[P\{n\}]= q^n [P],$ where $\{n\}$ is the grading shift by $n$ degrees up. 

$K_0(H_m)$ is a free $\Z[q,q^{-1}]$-module on one generator $[P_{(m)}]$, where
$P_{(m)}:= H_m e_m \{\frac{m(1-m)}{2} \}$ is an indecomposable graded 
projective $H_m$-module, unique up to grading shifts and isomorphisms. 
Placing diagrams next to each other 
gives inclusions $H_n \otimes H_m \subset H_{n+m}$,  
which lead to induction and restriction functors between categories of graded  
$H_n\otimes H_m$-modules and $H_{n+m}$-modules. These functors preserve 
subcategories of finitely-generated projective modules and give us maps 
\begin{eqnarray*} 
   M &  : &  K_0(H_n)\otimes K_0(H_m) \lra K_0(H_{n+m}), \\
 \Delta  &  : &    K_0(H_{n+m}) \lra \ K_0(H_n)\otimes K_0(H_m), 
\end{eqnarray*} 
where the tensor products are over $\Z[q,q^{-1}]$. Note that 
\begin{equation} \label{k0-iso} 
K_0(H_n \otimes_{\Bbbk} H_m) \ \cong \ K_0(H_n) \otimes_{\Z[q,q^{-1}]} 
 K_0(H_m)  , 
\end{equation} 
essentially due to absolute irreducibility of graded simple modules over these rings 
for any field $\Bbbk$, allowing us to freely switch between the two sides 
of (\ref{k0-iso}) and to define $\Delta$ (a similar fact for rings $R(\nu)$ was 
glossed over in~\cite{KL1}).  

Summing over all $n$ and $m$ produces maps that turn 
$$ K_0(H) \ := \ \oplusop{m\ge 0} K_0(H_m) $$ 
into a twisted $\Z[q,q^{-1}]$-bialgebra 
naturally isomorphic to $U^+_{\Z}(\mf{sl}(2)).$ 
This isomorphism takes $[H_m]$ to $E_2^m$ and $[P_{(m)}]$ to the 
divided power $E_2^{(m)}$. Following our notations we think of 
$U^+_{\Z}(\mf{sl}(2))$ as a $\Z[q,q^{-1}]$-subalgebra of
$U^+_{\Z}=U^+_{\Z}(\mf{gl}(1|2))$ generated by the divided powers of $E_2$. 

\vspace{0.2in} 


\noindent 
{\bf Lipshitz-Ozsv\'ath-Thurston dg algebras} 


Continuing to work over a field $\Bbbk$, consider the $\Bbbk$-algebra $\LOT_n$ 
with generators $\sigma_1, \dots, \sigma_{n-1}$ and defining relations 
\begin{eqnarray} 
   \sigma_i^2 & = & 0 ,  \label{eq-odd1}\\
  \sigma_i \sigma_j + \sigma_j \sigma_i & = & 0  \quad \quad \mbox{if} 
\quad \quad |i-j|>1, \label{eq-odd2}\\
    \sigma_i \sigma_{i+1}\sigma_i & = & \sigma_{i+1}\sigma_i \sigma_{i+1}. 
   \label{eq-odd3}
\end{eqnarray} 
Algebras $\LOT_n$ can be given a graphical description~\cite{LOT1}, by considering 
diagrams of $n$ lines 
in the strip $\R\times [0,1]$ of the plane, connecting $n$ fixed points on the bottom 
line $\R\times \{0\}$ to $n$ matching points on the top line $\R\times \{1\}$. 
Each line projects bijectively onto $[0,1]$ along $\R$. 
Diagrams are composed via concatenation, and isotopies are allowed that don't 
change the relative height position of crossings. In addition, the following relations 
are imposed 

\begin{equation} \label{eq-draw} 
\epsfig{file=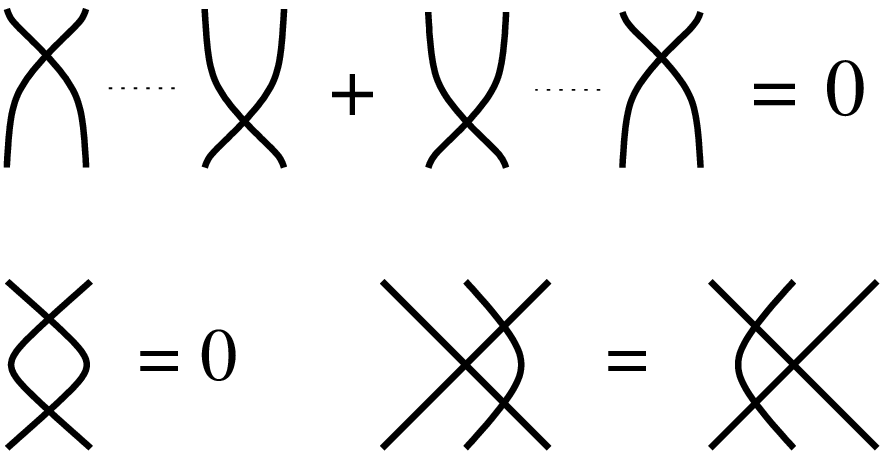}
\end{equation}

We draw these as thick lines, to distinguish them from the lines that enter 
the diagrammatics for the nilHecke algebra $H_m$. 
One can think of these lines as being fermionic, so that the far away crossing points 
anticommute rather than commute (change in the order of the product 
$\sigma_i \sigma_j\longleftrightarrow \sigma_j \sigma_i$ 
corresponds to the relative height change of the two far away crossings). 

We make $\LOT_n$ graded, with $\deg(\sigma_i)=-1$, thinking about it as cohomological 
grading (as opposed to the grading of $H_m$ above, which we'll refer to as $q$-grading), 
and equip $\LOT_n$ with a differential $d$ via the rules
$$ d(\sigma_i) = 1, \quad \quad d(ab) = d(a)b + (-1)^{\deg(a)} a d(b).$$ 
The differential takes a diagram to the alternating sum of diagrams obtained by 
smoothing a crossing $c$ of the diagram 
and multiplying what's left by $-1$ to the power the number of crossings above $c$.  

\drawing{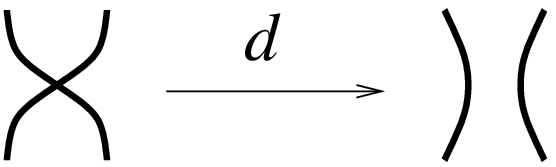}

This turns $\LOT_n$ into a differential graded (dg) $\Bbbk$-algebra. 
Ring $\LOT_n$ is the simplest example in the family of rings introduced by 
Lipshitz, Ozsv\'ath, and Thurston~\cite[Section 3]{LOT1}, \cite{LOT2}  
to extend Ozsv\'ath-Szab\'o 3-manifold homology to 3-manifolds with boundary 
and to localize combinatorial (grid diagram) construction~\cite{MOS}  
of the Ozsv\'ath-Szab\'o-Rasmussen knot Floer homology (which 
categorifies the Alexander polynomial). The authors of~\cite{LOT1, LOT2} work over 
$\Z/2$; the above characteristic-free lifting of their dg ring is a straightforward guess. 

\vspace{0.06in} 

It is obvious that the dimension of $H_n^-$ is at most $n!$ and slightly 
less obvious that the dimension is exactly $n!$. Let us explain this fact. 
A permutation $w\in S_n$ admits (usually non-unique) reduced expression 
$w= s_{i_1}\dots s_{i_r}$ as a product of transpositions $s_i = (i,i+1)$, 
with $r= {\it l}(w)$ the length of $w$, also equal to the number of crossings 
in any minimal presentation of $w$ via $n$ intersecting lines in the plane.  
We can describe a presentation $w'$ of 
$w$ by listing the sequence of indices  $w'=(i_1, \dots , i_r)$. To each presentation 
$w'$ we assign the element $\sigma_{w'} = \sigma_{i_1}\dots \sigma_{i_r}$ of 
$H_n^-$. Fixing a presentation for each permutation $w$, we obtain a set 
of elements $\{ \sigma_{w'}\}_{w\in S_n}$ that clearly spans $H_n^-$.

\begin{lemma} This set is a basis of $H_n^-$ as a $\Bbbk$-vector space. 
\end{lemma} 

\noindent 
\emph{Proof of lemma:} The only potential issue is that minus signs in 
the relations $\sigma_i \sigma_j = - \sigma_j \sigma_i$ for $|j-i|>1$ 
might force the relation $\sigma_{w'}=-\sigma_{w'}$ for some permutation 
$w$, making $\sigma_{w'}=0$ if $\mathrm{char}(\Bbbk)\not= 2$.  
To see that this does not happen, denote by 
$PDI(w)$ the set of \emph{pairs of disjoint inversions} in $w$. An inversion in $w$ is a pair 
$(j_1,j_2)$ of numbers such that $j_1<j_2$ but $w(j_1)>w(j_2)$. 
Given a reduced presentation $w'$ of $w$,  
inversions are in a bijection with terms of the presentation.  
If $w'$ is drawn as a diagram of $n$ intersecting lines in the plane, 
an inversion corresponds to a pair of intersecting lines (a crossing). 
A pair of disjoint inversions 
is a quadruple $(j_1,j_2,k_1, k_2)$ such that $(j_1,j_2)$ and $(k_1,k_2)$ are 
inversions, $j_1<k_1,$ and all four numbers $j_1,j_2,k_1,k_2$ are distinct. 
Diagrammatically, a pair of inversions corresponds to a pair of crossings 
that belong to four distinct lines. 

\drawing{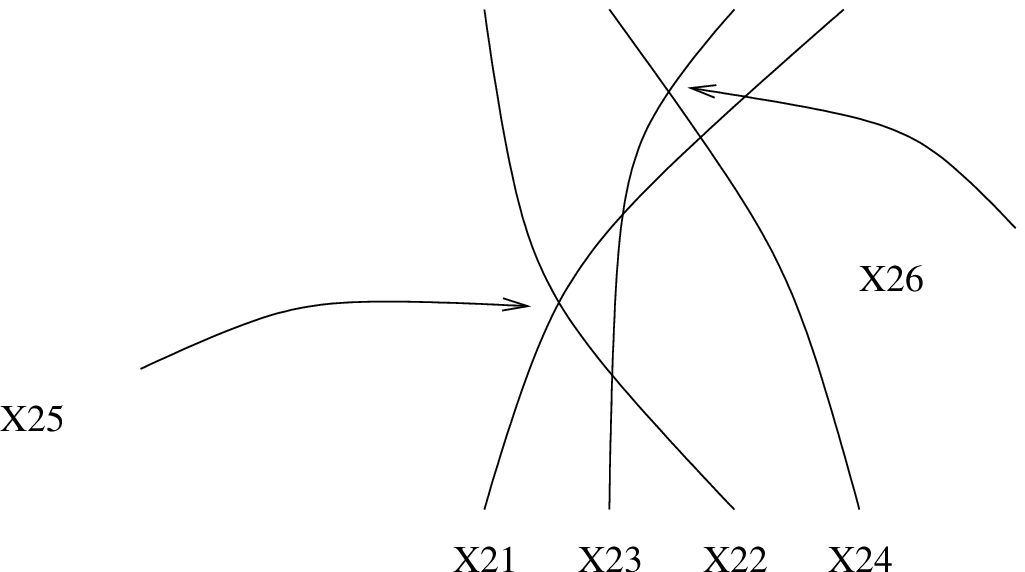}
 
Given $x=(j_1,j_2,k_1,k_2)\in PDI(w)$ and a presentation $w'$ of $w$, define 
${\epsilon}(w',x) =1 $  if the $(j_1,j_2)$-crossing is located 
below the $(k_1,k_2)$-crossing in the diagram of $w'$ and 
${\epsilon}(w',x) =-1$ if it is located above the $(k_1,k_2)$-crossing. 
Algebraically, ${\epsilon}(w',x) =1 $ if the generator for the 
$(j_1,j_2)$-crossing appears in $w'$ to the right of the generator 
for the $(k_1,k_2)$-crossing and  ${\epsilon}(w',x) =-1$ if it 
appears to the left. In the above diagram, ${\epsilon}(w',x) =1$. 

Define ${\epsilon}(w')\in \{ 1,-1\}$ as the product of ${\epsilon}(w',x)$ 
over all $x\in PDI(w)$, and let 
$$\widetilde{\sigma}_{w'} = {\epsilon}(w') \sigma_{w'}.$$ 
If $w''$ differs from $w'$ by a simple transposition of two consecutive terms,  
$$ w'=(\dots, i,j,\dots ),\quad \quad w''=(\dots, j, i, \dots ),\quad \quad |j-i|>1$$ 
(geometrically, $w',w''$ differ as the two diagrams in the top equation in 
(\ref{eq-draw})), then 
${\epsilon}(w',x) = {\epsilon}(w'',x)$ for all $x\in PDI(w)$ save the  
one that corresponds to the permuted pair of crossings. Hence, 
${\epsilon}(w'') = - {\epsilon}(w')$, matching the sign change in the 
equation $\sigma_j \sigma_i = - \sigma_i \sigma_j$,  and 
$\widetilde{\sigma}_{w'} =\widetilde{\sigma}_{w''}$
in this case. 

If $w''$ differs from $w'$ by a "Reidemeister III" move, 
$$ w'=(\dots, i,i+1,i,\dots ),\quad \quad w''=(\dots, i+1, i, i+1,\dots ),$$ 
then ${\epsilon}(w',x) = {\epsilon}(w'',x)$ for all $x\in PDI(w)$, 
and $\widetilde{\sigma}_{w'} = \widetilde{\sigma}_{w''} $. 

These observations imply that if $w'$ and $w''$ are two presentations of $w$  then 
 $\widetilde{\sigma}_{w'} = \widetilde{\sigma}_{w''}$, so that we can define 
${\sigma}_w := \widetilde{\sigma}_{w'} $ for any presentation $w'$ of 
$w$, and $\sigma_w$ will depend only of $w$ and not on its presentation. 
 That ${\sigma}_w\not= 0$ (consistency) follows as well. $\square$ 

\vspace{0.06in} 

\emph{Remark:} If $\mathrm{char}(\Bbbk)=2$, anticommutativity is 
indistinguishable from commutativity, and $\LOT_n$ turns into 
the nilCoxeter algebra, a subalgebra of the nilHecke algebra generated by the divided 
difference operators. The nilCoxeter algebras (over any field and without the differential) 
can be used to categorify the polynomial representation of the first Weyl algebra, 
as well as the bialgebra $\Z[E], \Delta(E)=E\otimes 1 + 1\otimes E$, see~\cite{Kh1}
(as opposed to the bigger bialgebra $\Z[E^{(m)}]_{m\ge 1}$ that contains divided powers 
$E^{(m)}=\frac{E^m}{m!}$ of $E$ and whose 
categorification relies on nilHecke rings $H_m$, see earlier but without the $q$-grading). 

\vspace{0.06in} 

Placing diagrams next to each other gives dg ring inclusions 
\begin{equation}
\LOT_n \otimes \LOT_m \subset \LOT_{n+m} \label{eq-inclusions} 
\end{equation} 
(the tensor product is taken in the category of super vector spaces, and a 
diagram from $\LOT_n$ is placed above and to the left of a diagram from $\LOT_m$). 
We would like to form the Grothendieck group $K_0(\LOT_n)$ and use these inclusions 
to define a multiplication and comultiplication on 
$$ K_0(\LOT) \ := \ \oplusop{n \ge 0} K_0( \LOT_n), $$ 
then identify $K_0(\LOT)$ with the integral subalgebra of the positive half of 
quantum $\mf{gl}(1|1)$. This integral subalgebra is $\Z[q,q^{-1},E_1]/(E_1^2)$.  

\vspace{0.2in} 


\noindent 
{\bf $K_0$ of a dg ring} 


\noindent 
The analogue of a projective module over a $\Bbbk$-algebra $A$ is 
a projective dg module over a dg $\Bbbk$-algebra (called $\mc{K}$-projective in~\cite{BL}). 
A (left) dg module $M$ over $A$ is a $\Z$-graded $A$-module equipped with a 
differential $d_M : M^i \lra M^{i+1}$ such that $d_M(am)= d(a)m + (-1)^{\deg(a)} 
a d_M(m)$. We call a dg module $P$ over a dg $\Bbbk$-algebra $A$  \emph{projective} if 
the complex $\Hom_A(P,M)$ has 
zero homology whenever $M$ does; here $M$ is a dg module over $A$. 

For an introduction to dg modules and projective dg modules we refer the 
reader to~\cite[Section 10]{BL}. Starting with the abelian category of dg $A$-modules, 
one first produces a triangulated category $\mc{K}(A)$ of dg-modules by 
modding out by homotopic to zero morphisms, then quasi-isomorphisms are inverted to 
produce the derived category $D(A)$.  The category $\mc{K}(A)$ contains 
the full subcategory $\mc{KP}(A)$ of projective dg modules. The localization 
functor $\mc{K}(A) \lra D(A)$, when restricted to $\mc{KP}(A)$, gives 
an equivalence $\mc{KP}(A)\cong D(A)$. 

To define the Grothendieck group $K_0(A)$ we need to restrict the size of 
projective modules. An object $M$ of $\mc{KP}(A)$ or $D(A)$ is called \emph{compact} 
if the inclusion 
$$ \oplusop{i\in I} \Hom (M, N_i) \subset \Hom(M, \oplusop{i\in I}N_i) $$ 
is an isomorphism for any collection $\{N_i\}_{i\in I}$ of objects indexed by 
a set $I$. This definition of a compact object makes sense in any additive 
category which admits arbitrary direct sums (not just finite ones). 
In a category of modules over a ring $A$, a module is compact iff it is 
finitely-generated as an $A$-module. 

Let $\mc{P}(A) \subset \mc{KP}(A)$ be the full subcategory 
of compact projective modules, for a dg algebra $A$. It is a triangulated 
category. Define $K_0(A)$  as the Grothendieck group 
of $\mc{P}(A)$. 
It has generators $[P]$ over all compact projectives 
$P$ and relations $[P[1]]= -[P]$ (here $[1]$ is the grading shift), and 
$[P_2] = [P_1]+[P_3]$ for each distinguished triangle $P_1 \lra P_2 \lra P_3$. 
Note that $\mc{P}(A)$ is equivalent to the subcategory $\mc{P}'(A)$ of 
compact objects in $D(A)$, and 
we can alternatively define $K_0(A)$ as the Grothendieck group of 
the triangulated category $\mc{P}'(A)$.  The following diagram summarizes 
our categories, inclusions, and equivalences. 

\begin{equation} \label{diag-cd} 
\xymatrix{
 A\mathrm{-dgmod}  \ar[r]   & \mc{K}(A) \ar[rd]  &    \\
   &     \mc{KP}(A) \ar@{^{(}->}[u] \ar[r]^{\cong} & D(A) \\
   &    \mc{P}(A)   \ar@{^{(}->}[u]  \ar[r]^{\cong} & \mc{P}'(A) \ar@{^{(}->}[u] 
} 
\end{equation}

A quasi-isomorphism $A\lra B$ of dg algebras induces an equivalence 
of derived categories $D(A)\cong D(B)$, an equivalence of subcategories 
of compact objects $\mc{P}'(A) \cong \mc{P}'(B)$, and an isomorphism 
$K_0(A)\cong K_0(B)$.

Any $\Bbbk$-algebra $A$ is naturally a dg algebra, concentrated in degree $0$ and with the 
trivial differential. In this case, in addition to the above definition of $K_0(A)$, 
there is the classical definition of $K_0(A)$ as the Grothendieck group of
projective finitely-generated $A$-modules. 

\begin{lemma} If $A$ is (left) Noetherian, the two definitions give naturally 
isomorphic groups. 
\end{lemma} 

The lemma can be proved by showing that $\mc{P}'(A)$ is equivalent to the 
homotopy category of bounded complexes of finitely-generated projective left 
modules over the ring $A$. $\square$ 

For later use, we point out that the above story has a generalization if 
the dg algebra $A$ has an additional grading, complementary to the cohomological 
grading. We call such $A$ a graded dg ring and refer to the additional grading 
as $q$-grading. Then one can talk about the category of graded dg modules, its 
homotopy and derived categories, projective graded dg modules, etc. 
Retaining the above notations, the Grothendieck group $K_0(A)$ will be
 a $\Z[q,q^{-1}]$-module, with $q$ corresponding to the grading shift in 
the additional grading. If $A$ is just a graded algebra, we turn it into a graded dg algebra 
by placing it entirely in cohomological degree $0$, so that the differential 
acts by $0$, and form the $\Z[q,q^{-1}]$-module $K_0(A)$, the 
Grothendieck group of the category of compact objects in $D(A)$, the derived 
category of the category of graded dg $A$-modules. 
The classical definition of $K_0$ of a graded ring, mentioned earlier (in our 
discussion of $H_m$), 
also produces a $\Z[q,q^{-1}]$-module, with generators $[P]$, over 
graded finitely-generated projective $A$-modules, and relations coming from 
direct sum decompositions. 

\begin{lemma} If $A$ is (left) graded Noetherian, the two definitions give naturally 
isomorphic $\Z[q,q^{-1}]$-modules. \label{lemma-noeth}   
\end{lemma} 

We say that a (graded) dg algebra $A$ is (graded) formal if it is (graded) quasi-isomorphic 
to its cohomology algebra $H(A)$.  In this case we can 
identify $K_0(A) \cong K_0(H(A))$. An easy exercise shows that $A$ is formal if 
$H(A)$ is concentrated in cohomological degree $0$. If, furthermore, $H(A)$ is a 
(graded) Noetherian algebra, we can describe $K_0(A)$ via finitely-generated 
(graded) projective $H(A)$-modules. 

From now on all algebras and dg algebras that we consider are graded, and we 
work with the category of graded dg modules, its homotopy and derived 
categories. Corresponding $K_0$-groups are $\Z[q,q^{-1}]$-modules. 
The $q$-grading on $\LOT_n$ is trivial -- the entire dg algebra sits in zero $q$-degree. 

\vspace{0.2in} 


\noindent 
{\bf Categorification of positive half of quantum $\mf{gl}(1|1)$} 


\noindent 
Let us compute $K_0$ of rings $\LOT_n$. The ring $\LOT_0=\Bbbk$
(the only diagram when $n=0$ is the empty one), and
$K_0(\LOT_0) \cong \Z[q,q^{-1}]$, with the generator $[\Bbbk]$, 
since an object in the category of complexes of graded vector spaces up to 
chain homotopy is compact iff its total cohomology is finite-dimensional.  
The ring $\LOT_1= \Bbbk$, since when $n=1$ the diagrams have only 
one line and no room for interactions. Again, $K_0(\LOT_0) \cong \Z[q,q^{-1}]$. 

To treat the $n\ge 2$ case we use the following observation. 

\begin{lemma} \label{throw-idemp} 
Suppose that $A$ is a dg ring and $x\in A$ an element of degree $-1$ such that 
$e=dx$ is an idempotent. Then the inclusion $ (1-e) A(1-e) \subset A$ 
is a quasi-isomorphism. In particular, if $dx=1$ for some $x\in A$ then 
$\mathrm{H}(A)=0$. 
\end{lemma}
\noindent 
\emph{Proof:} 
An idempotent $e$ in a ring $A$ makes it look (superficially) like the 
ring of $2\times 2$-matrices: 
$$ A = eAe\oplus eA(1-e) \oplus (1-e)Ae\oplus (1-e)A(1-e) .$$ 
For $e$ and $A$ as in the lemma, each of the first three summands  
is a contractible complex of abelian groups. The map 
$h: Ae \lra Ae$ that takes $ae$ to $(-1)^{|a|}a x e$ satisfies $hd+ dh=1$, 
implying that $Ae$ is contractible; a similar computation establishes contractibility 
of the second summand. $\square$ 

\vspace{0.05in} 

We call idempotents that can be written as $dx$, for some $x$, \emph{contractible} 
idempotents. 

\vspace{0.1in} 

The dg ring $\LOT_n$ has a special property that, for $n\ge 2$, the differential of some 
element equals one, for instance $d(\sigma_i)=1$. This implies that $\LOT_n$ has trivial 
homology and its derived category $D(\LOT_n)$ is trivial. Category 
$\mc{P}(A)$ is equivalent to $\mc{P}'(A)$, a full subcategory of $D(\LOT_n),$ 
so it is trivial as well and $K_0(\LOT_n)=0$ for $n\ge 2$. 

This way we do obtain  $\Z[q,q^{-1},E_1]/(E_1^2)$, the integral form of 
$U^+_q(\mf{gl}(1|1))$, as the Grothendieck group 
$$K_0(\LOT)\ := \ \oplusop{n\ge 0} K_0(\LOT_n).$$ 
Generators $1, E_1\in \Z[q,q^{-1},E_1]/(E_1^2)$ are given by the symbols $[\LOT_0]$ 
and  $[\LOT_1]$ of free modules over $\LOT_0\cong \Bbbk$, $\LOT_1\cong \Bbbk$, 
respectively, and multiplication, comultiplication descend from the induction and restriction 
functors. 

Notice that we went into a lot of trouble, only to categorify 
the exterior algebra on one generator, with a bit of additional structure thrown in. 
This will pay off momentarily, as well as teach us that to categorify algebras with 
nilpotent generators it helps to use categories of modules over dg rings rather than 
just of modules over rings. 

\vspace{0.15in} 


\noindent 
{\bf Putting the two categorifications together.} 


\noindent 
The Lie superalgebra $\mf{gl}(1|2)$ is generated by its subalgebras $\mf{gl}(1|1)$ 
and $\mf{gl}(2)$. Likewise, $U^+_q=U^+_q(\mf{gl}(1|2))$ has generators $E_1, E_2$
(see the first section), first odd, second even, which generate a copy of 
$U^+_q(\mf{gl}(1|1)),$ $U^+_q(\mf{sl}(2))$, respectively. We know categorifications 
of these subalgebras, via diagrammatics of braid-like pictures, and simply need to 
guess how to combine them. A possible answer is the following. 

For each $n,m\in \Z_+$ let $R(n,m)$ be the dg $\Bbbk$-algebra spanned by 
braid-like diagrams with $n$ lines labelled 1 and $m$ lines labelled 2. 
We call lines of the first type \emph{fermionic} or \emph{odd}, lines of the second type 
\emph{bosonic} or \emph{even}. 
We draw fermionic lines thicker than bosonic. Bosonic lines can carry dots, fermionic 
lines cannot. Bosonic lines and dots on them interact via the diagrammatics for 
the nilHecke algebra, fermionic lines interact through diagrammatics for the LOT 
(Lipshitz-Ozsv\'ath-Thurston) dg algebra. 
Alternatively, we will indicate fermionic lines by labelling their lower endpoints $1$ 
and bosonic lines by labelling their lower endpoints $2$ (as in generators $E_1$, $E_2$
of $U^+_q$). The following example displays the relation between the two notations: 

\drawing{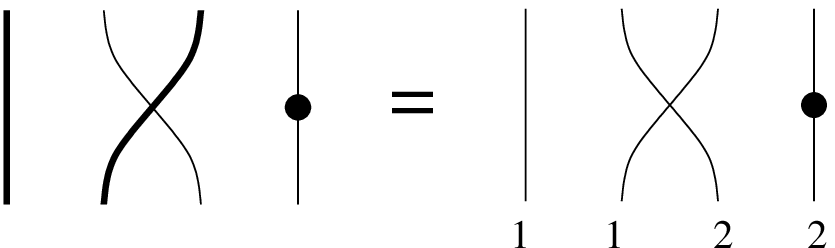}

\noindent 
We need to add additional generators -- intersections between 
fermionic and bosonic lines: 

\drawing{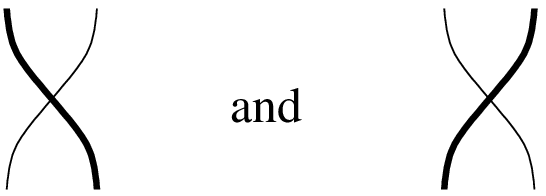}

\noindent 
and impose the following relations: 

1) Far away intersections commute, unless both intersections are between fermionic 
lines, in which case they anticommute. We can encode these into a single relation 

\drawing{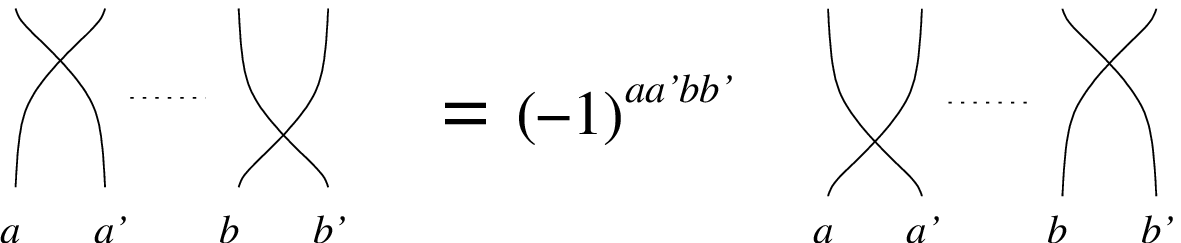}

\noindent 
where $a, a', b, b'\in \{1,2\}$ are the labels of the four crossing lines. 
The only anticommuting case is the following: 

\drawing{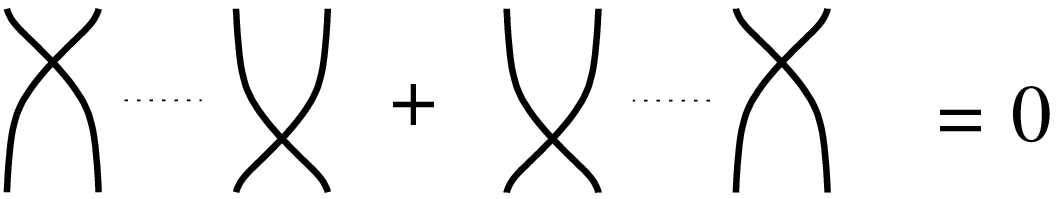}

2) A dot commutes with far away intersection of any kind

\drawing{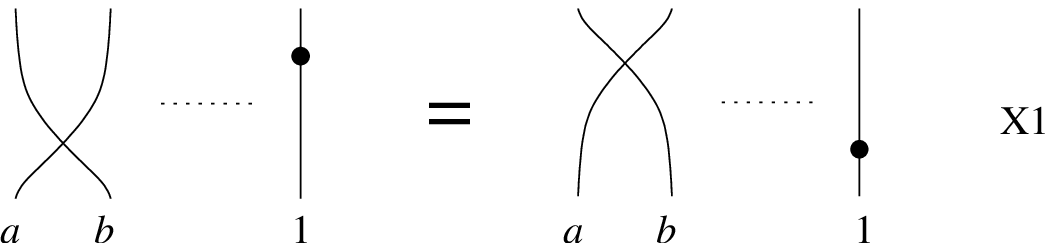}

\noindent 
(same if the intersection is to the right of the dot). Two dots commute: 

\drawing{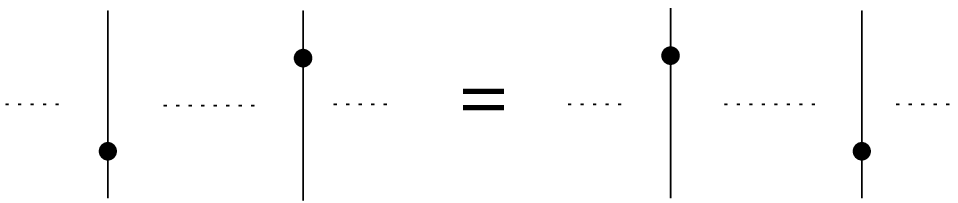} 

\noindent 
Summary of relations 1) and 2): far away generators commute if at 
least one of them is even and anticommute if both are odd. 

\vspace{0.06in} 

3) Dot through a crossing relations: 

\drawing{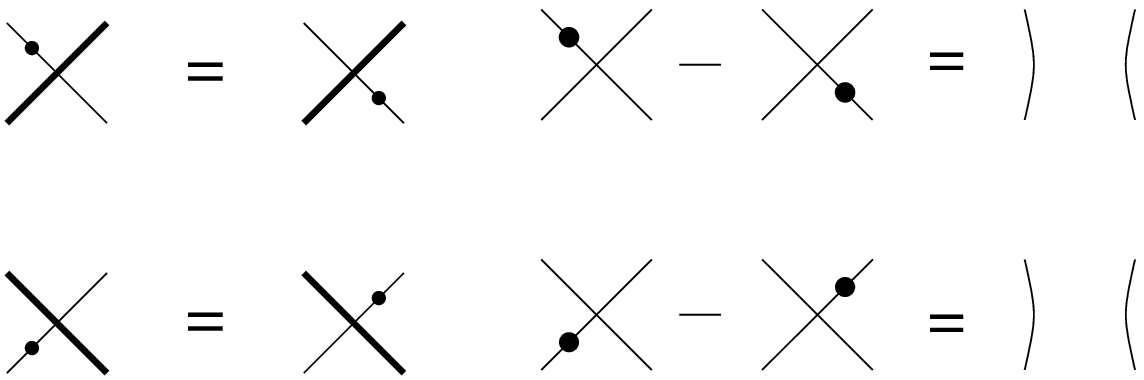}

\noindent 
These say that a dot can freely move through an odd-even crossing, 
and can move through an even-even crossing at the cost of adding an extra term, 
as in the nilHecke algebra. 

\vspace{0.06in} 

4) Two-line relations (Reidemeister II type relations): 

\drawing{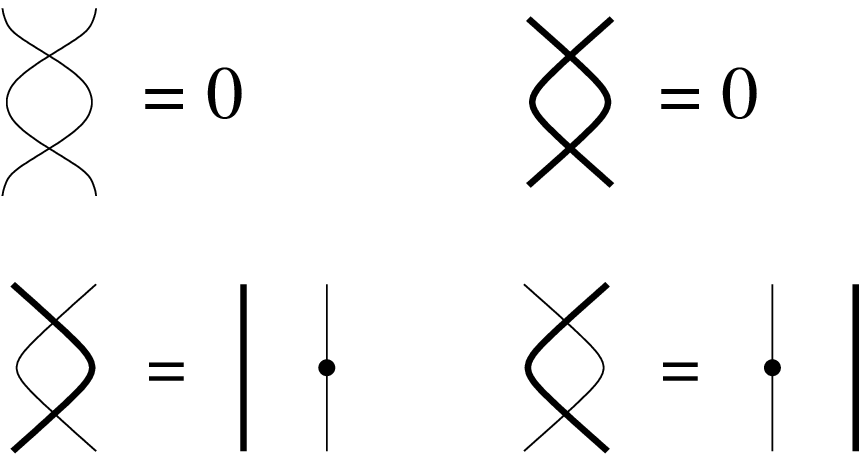}

\noindent 
If the two lines that make a crossing are both bosonic or both fermionic, the 
square of the crossing is $0$. A double crossing of a bosonic and a fermionic 
line equals to the vertical lines diagram with a dot on the bosonic line. 

\vspace{0.06in} 

5) Three-line relations (Reidemeister III type relations):  

\drawing{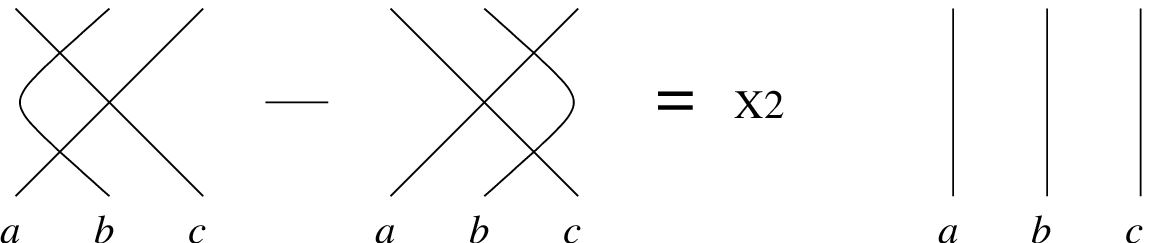}

\noindent 
These relations say that a triple intersection homotopy is allowed, unless 
the three line types are even, odd, even (in this order), in which case 
there is an additional term:  

\drawing{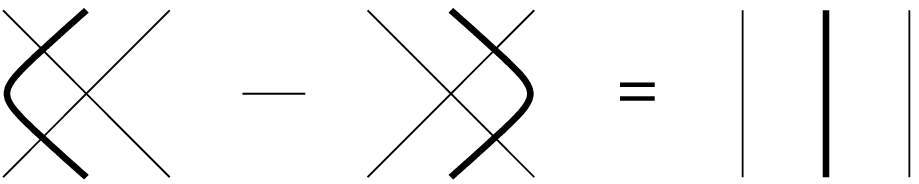}

The defining relations for $R(n,m)$ contain the nilHecke and the LOT relations. 
When bosonic lines are absent (case $m=0$), they are exactly the 
relations in the LOT algebra, so that $R(n,0)\cong \LOT_n$.
When fermionic lines are absent ($n=0$ case), the relations are that of 
the nilHecke algebra, and $R(0,m)\cong H_m$. 

For each sequence ${\bf i}$ of $n$ ones and $m$ twos we have the corresponding 
idempotent $1_{\bf i}$ in $R(n,m)$ with the diagram of $n+m$ vertical lines, $n$ fermionic and 
$m$ bosonic, in the order ${\bf i}$. The unit element of $R(n,m)$ is the sum of 
these idempotents, over all sequences ${\bf i}$. Let $P_{\bf i}= R(n,m)1_{\bf i}$ 
be the left projective $R(n,m)$-module corresponding to the idempotent $1_{\bf i}$.  

We make $R(n,m)$ bigraded by placing each $1_{\bf i}$ in bidegree $(0,0)$ 
and listing bidegrees of other generators in the table below: 

\begin{equation}  \label{eq-table} 
\epsfig{file=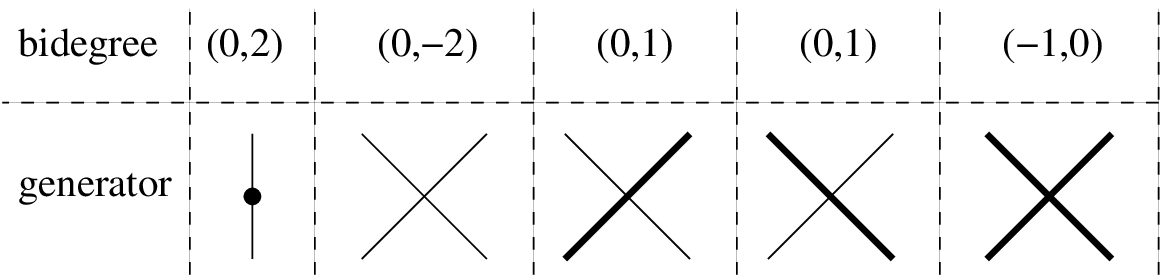}
\end{equation} 

Defining relations are homogeneous, and $R(n,m)$ is a bigraded $\Bbbk$-algebra. 
We call the first grading \emph{cohomological grading}, the second grading \emph{$q$-grading}. 
This bigrading restricts to the previously 
discussed gradings on the nilHecke and LOT algebras.  

We turn $R(n,m)$ into a dg algebra by defining $d$ on generators (crossings, dots, 
idempotents $1_{\bf i}$) to be $0$, except 
that $d$ of an odd crossing is the idempotent $1_{\bf i}$ given by resolving the crossing:  

\drawing{differential.eps}

\noindent 
Bidegree of the differential is $(1,0)$, thus $d$ respects the $q$-grading, 
and $R(n,m)$ becomes a graded dg algebra.  

\vspace{0.06in} 

\noindent 
Putting diagrams next to each other defines inclusions of graded dg algebras 
$$ R(n,m)\otimes R(n',m') \subset R(n+n', m+m') $$
(with the same caveats as for inclusions (\ref{eq-inclusions})) 
and leads to induction and restriction functors between categories of graded 
dg modules, and between corresponding derived categories. We claim that 
$K_0(R(n,m))$ can be identified with the weight $n\alpha_1 + m\alpha_2$ 
subspace of $U^+_{\Z}$, 
\begin{equation} \label{K0-eqn} 
K_0(R(n,m)) \cong U^+_{\Z} (n,m), 
\end{equation} 
so that induction and restriction functors descend to multiplication and 
comultiplication on $U^+_{\Z}$: 
\begin{eqnarray*} 
K_0(R) & \cong & U^+_{\Z}, \\
K_0(R)  & := &  \oplusop{n,m\ge 0}K_0(R(n,m)), \quad \quad 
   U^+_{\Z} = \oplusop{n,m\ge 0}U^+_{\Z}(n,m). 
\end{eqnarray*} 
(do not confuse the notation $(n,m)$ for the weight space 
of $\mf{gl}(1|2)$ with the earlier notation $(n|m)$ for the parameters 
of $\mf{gl}(n|m)$). 

$P_{\bf i}=R(n,m)1_{\bf i}$ is a dg module, since $d(1_{\bf i})=0$. We want 
isomorphism (\ref{K0-eqn}) to take 
$[P_{\bf i}]$  to $E_{\bf i}$, were $E_{\bf i} = E_{i_1}\dots E_{i_{n+m}}$ for 
${\bf i}=i_1i_2\dots i_{n+m}$. 
Also, the divided power element $E_2^{(m)}$ should correspond  
to the indecomposable projective $R(0,m)$-module 
\begin{equation} \label{eq-p2m}
P_{2^{(m)}} \ := \   R(0,m)e_m \left\{ \frac{m(1-m)}{2} \right\}.
\end{equation}  
The idempotent $e_m\in H_m$ was described earlier. Here we use the same 
notation $e_m$ for the image of $e_m$ in $R(0,m)$ under the 
canonical isomorphism $H_m \cong R(0,m)$.  

\vspace{0.06in} 

Reflecting diagrams about the $x$-axis induces an anti-involution on $R(n,m)$. 
Reflecting a diagram about the $y$-axis and multiplying it by $(-1)^b$, where 
$b$ is the number of 22-crossings, induces an involution on $R(n,m)$. 

\vspace{0.06in} 

Any diagram representing an element of $R(n,m)$ can be simplified to a linear combination 
of diagrams which consist of fixed minimal length presentations of permutations $w\in S_{n+m}$ 
with some number of dots at the top of each bosonic strand. Similar to the case of 
rings $R(\nu)$, see~\cite{KL1}, this set is a basis of $R(n,m)$ as a $\Bbbk$-vector space. 
For rings $R(\nu)$ this was shown~\cite{KL1} 
by checking that these elements acts linearly independently 
on a certain representation $Pol_{\nu}$ of $R(\nu)$. 

\vspace{0.06in} 

The analogue of $Pol_{\nu}$ is a representation $Pol$ of $R(n,m)$ given 
by 
$$Pol= \oplusop{\ii, w} Pol(\ii,w), \quad \quad 
 Pol(\ii,w) = \Bbbk[x_1(\ii,w),   \dots , x_m(\ii,w)].$$ 
Here $\ii$ ranges over all sequences of $n$ ones and $m$ twos, and 
$w$ over elements of the symmetric group $S_n$. The idempotent $1_\ii$ acts 
as identity on $Pol(\ii,w)$ for each $w$ and by $0$ on  
$Pol(\jj ,w')$ for $\jj \not= \ii$. An element $x\in  R(n,m)1_{\ii}$ takes 
elements of $Pol(\jj,w)$ to $0$ for $\jj\not= \ii$. 

A diagram of vertical lines for a sequence $\ii$ with a dot on the $i$-th bosonic 
strand counting from the left takes $f\in Pol(\ii,w)$ to $x_i(\ii, w)f $
(these diagrams are generators of the first type listed in table (\ref{eq-table})). 

A crossing diagram of the $i$-th and $(i+1)$-st bosonic strands counting from the 
left (assuming they are next to each other) acts as the divided difference operator, taking 
$f\in Pol(\ii,w)$ to $\partial_i(f)\in Pol(\ii, w)$. These diagrams are generators 
of the second type listed in table (\ref{eq-table}). 

A crossing diagram of the $i$-th bosonic and $j$-th fermionic strands, assuming 
they are next to each other in the sequence and 
the bosonic strand is on the main diagonal (third type generator in the table) takes 
$f \in Pol(\ii, w)$ to the same polynomial $f$ in variables $x_1(\ii', w), \dots, 
x_m(\ii',w)$ instead of $x_1(\ii, w)$,  $\dots, x_m(\ii,w)$, where $\ii'$ is given by 
transposing $i$-th two with $j$-th one in the sequence $\ii$: 
$$ \ii = \dots 12 \dots \quad , \quad \quad \ii'= \dots 21 \dots \quad .$$  

A crossing diagram of the $i$-th bosonic and $j$-th fermionic strands, assuming 
they are next to each other in the sequence and 
the fermionic strand is the main diagonal (fourth type generator in the table) takes 
$f \in Pol(\ii, w)$ to $x_i(\ii',w) f'$, where $f'$ is the same as $f$ but with 
variables $x_1(\ii', w), \dots, 
x_m(\ii',w)$ substituted for  $x_1(\ii, w), \dots, x_m(\ii,w)$. Here $\ii'$ is given by 
transposing $i$-th two with $j$-th one in the sequence $\ii$: 
$$ \ii = \dots 21 \dots, \quad \quad \ii'= \dots 12 \dots .$$  

The crossing diagram of the $i$-th and $(i+1)$-st fermionic strands (assuming 
they are next to each other; this is a fifth type generator from (\ref{eq-table})) 
takes $f\in Pol(\ii,w)$ to 
$\epsilon_w^i f'\in Pol(\ii, s_iw)$ if ${\it l}(s_i w) = {\it l}(w) + 1$ (to go from $f$ to 
$f'$ change variables $x_i(\ii,w)$ to $x_i(\ii, s_i w)$) and to $0$ if 
${\it l}(s_i w) = {\it l}(w) - 1$, 
where ${\it l}$ is the usual length function in the symmetric group, and 
$\epsilon_w^i\in \{1, -1\}$ is determined by the formula 
$\sigma_{s_i w} = \epsilon_w^i\sigma_i \sigma_w$ (see the definition 
of $\sigma_w$ earlier). 

\vspace{0.06in} 

It is an easy but  enlightening exercise to check that these rules give an action 
of $R(n,m)$ on $Pol$. An argument similar to the one in \cite[Section 2.3]{KL1} 
shows that the action is faithful and implies that the spanning set in $R(n,m)$ described 
earlier is a basis of this algebra. 
If desired, the rules for the action of the two types of bosonic-fermionic 
crossings can be interchanged, resulting in the same endomorphism algebra. 

\vspace{0.06in} 

The algebra of symmetric polynomials $\mathrm{Sym}(n,m)$ in dots on bosonic lines 
belongs to the center of $R(n,m)$. Since $R(n,m)$ is a free module of finite rank over 
$\mathrm{Sym}(n,m)$, we conclude that $R(n,m)$ is both left and right Noetherian. 

\vspace{0.06in} 

We can construct a homomorphism 
\begin{equation} \label{eq-gamma} 
\gamma \ : \ U^+_{\Z} \lra K_0(R) 
\end{equation} 
by taking generator $E_1$ to $[P_1]$, where $P_1=R(1,0)$ is the 
free rank one $R(1,0)$-module, and generator $E_2^{(m)}$ to $[P_{2^{(m)}}]$. 
According to our definition, $U^+_{\Z}$ is spanned by 
various products $E_1 E_2^{(m_1)} E_1 \dots E_1 E_2^{(m_k)} E_1$ and their 
variations given by removing the first or the last $E_1$ from the product, or both. 
$\gamma$ will take the above product to the symbol of projective module 
$$P_{1 2^{(m_1)} 1 \dots 1 2^{(m_k)}1} \ := \ R(n,m) 1_{1 2^{(m_1)}1\dots 
1 2^{(m_k)}1 }  \{a \},$$ 
where $ 1_{1 2^{(m_1)}1\dots 1 2^{(m_k)}1 }$ is the idempotent obtained 
by placing  $e_{m_1}, \dots , e_{m_k}$ in parallel 
next to each other, separated by odd lines, one for each $E_1$ in the product, 
and $a= \sum_i \frac{m_i (1-m_i)}{2}$: 

\drawing{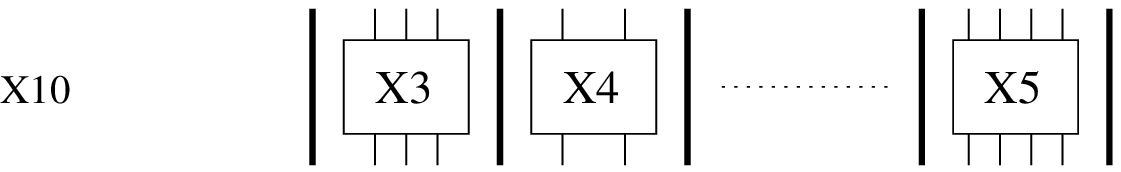} 

To check that $\gamma$ is well-defined  we consider the defining 
relations (\ref{eq-one}), (\ref{eq-two}). We already saw that $[P_{11}]=0$, 
due to contractibility of the dg module $P_{11}$. More generally, for any sequences 
$\ii, {\bf j}$ projective module $P_{\ii11{\bf j}}$ is contractible 
as a complex of $\Bbbk$-vector spaces, so that $[P_{\ii11{\bf j}}]=0$. 
In particular, the equivalent of relation (\ref{eq-one}) holds in $K_0(R)$. 

\vspace{0.15in} 


\noindent 
{\bf Comparison with categorified $\mf{sl}(3)$} 


\noindent 
To see that the equivalent of relation (\ref{eq-two}) holds in $K_0(R)$, 
we will compare $R(1,m)$ with the graded ring 
$R(\alpha_1 + m\alpha_2)$ that categorifies the weight space 
$\alpha_1 + m \alpha_2$ of $U^+_q(\mf{sl}(3))$. This is a
special case of rings $R(\nu)$ defined in~\cite{KL1}. The Dynkin diagram 
of $\mf{sl}(3)$ consists of two vertices joined by an edge, in this way similar 
to the Dynkin diagram of $\mf{gl}(1|2)$, which also consists of two vertices 
and an edge. Just like diagrams describing $R(1,m)$, diagrams for 
$R(\alpha_1 + m\alpha_2)$ have one line labelled $1$ 
and $m$ lines labelled $2$. But now the line labelled one can carry dots, 
which freely slide through intersections with type $2$ lines (there are no 
intersections of two lines labelled one, since in this weight there is only one 
such line). Another difference from $R(1,m)$ is that the double intersection 
of line 1 and line 2 equals the sum of two terms rather than one: a dot on line 
1 plus a dot on line 2:  

\drawing{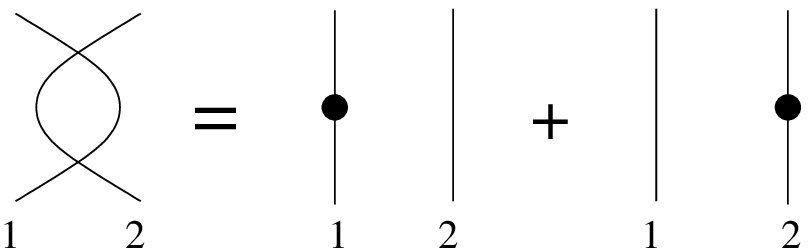} 

\noindent 
Other relations are the same. This gives a homomorphism 
$$\tau: R(\alpha_1 + m\alpha_2) \lra R(1,m) $$ 
that kills any diagram which contains a dot on type 1 line and is the identity on 
diagrams without such dots. The kernel of $\tau$ is spanned by diagrams 
with at least one dot on line 1. Equivalently, the kernel is the two-sided ideal 
generated by diagrams of vertical lines, with a dot on the line labelled 1.  

In the relations for $R(\alpha_1 + m\alpha_2)$ dot on line 1 can be slid 
up and down without obstacles. Let $J$ be the graded Jacobson radical 
of the graded ring $R(\alpha_1+m\alpha_2)$. We know from~\cite{KL1} 
that $J$ has finite codimension in $R(\alpha_1+m\alpha_2)$
and the quotient $R(\alpha_1+m\alpha_2)/J$ is a finite-dimensional 
semisimple $\Bbbk$-algebra. The 
ideal $(\mathrm{ker}(\tau))^N$ is spanned by diagrams with at least 
$N$ dots on line $1$, and its easy to see that, for degree reasons, 
 $(\mathrm{ker}(\tau))^N\subset J$ 
for sufficiently large $N$. This in turn implies that $\mathrm{ker}(\tau)\subset J$, 
so that the induced map of Grothendieck groups of graded rings  
$$K_0(\tau) \ : \ K_0(R(\alpha_1+m\alpha_2)) \lra K_0(R(1,m))$$
is an isomorphism. 

Given a sequence $\ii$ of divided powers of symbols $1$ and $2$, we 
associate to it an idempotent $1_\ii \in R(n\alpha_1 + m\alpha_2)$ and 
graded projective $R(n\alpha_1 + m\alpha_2)$-module
$$ P_\ii' = R(n\alpha_1 + m\alpha_2) 1_\ii\{a\} $$
for $n$ and $m$ equal to total weight of $1$ and $2$ in the sequence
(in~\cite{KL1} this module was  denoted $P_\ii$), with grading shift as 
in~\cite{KL1}. We only treat the case $n=1$, and then sequence $\ii$ also 
defines graded projective $R(1,m)$-module, denoted 
$$P_\ii\cong P_\ii'/\mathrm{ker}(\tau) P_\ii' = R(1,m) 1_\ii \{ a\},$$ 
where $a=\sum \frac{m_i(1-m_i)}{2}$ and $m_i$'s are divided powers of 2 that 
appear in $\ii$. 

Proposition 2.13 in~\cite{KL1} implies that, in our notations,  
$$ P'_{212} \cong P'_{1 2^{(2)}}\oplus P'_{2^{(2)} 1}.$$

Applying $\tau$, we see that 
$$ P_{212} \cong P_{1 2^{(2)}}\oplus P_{2^{(2)} 1},$$ 
and, more generally, 
$$ P_{\ii 212{\bf j}} \cong P_{{\bf i}1 2^{(2)}{\bf j}}\oplus 
P_{{\bf i}2^{(2)} 1{\bf j}},\quad \quad 
[ P_{{\bf i}212{\bf j}}]=[P_{{\bf i}1 2^{(2)}{\bf j}}] + 
  [P_{{\bf i}2^{(2)} 1{\bf j}}],$$ 
for any sequences ${\bf i}, {\bf j}$ of ones and twos. 
Therefore, relations (\ref{eq-two}), (\ref{eq-three}) hold in $K_0(R)$.

From~\cite[Section 3.3]{KL1} we know that $K_0(R(\alpha_1 + m\alpha_2))$
is a free $\Z[q,q^{-1}]$-module with basis elements $[P'_{1 2^{(m)}}]$, 
$[P'_{2^{(m)} 1 }]$ (this is implied by $\{ E_1 E_2^{(m)}, E_2^{(m)} E_1\}$ 
being the canonical basis of weight $\alpha_1+ m \alpha_2$ subspace of 
$U^+_q(\mf{sl}(3))$). Moreover, any finitely-generated graded projective 
$R(\alpha_1 + m\alpha_2)$-module is isomorphic to a direct sum of (graded shifts of) 
modules $P'_{1 2^{(m)}}, P'_{2^{(m)} 1 }$, with multiplicities 
determined by the image of the module in the Grothendieck group.  In particular,  
$$ P'_{2^{(k)}12^{(m-k)}} \cong (P'_{12^{(m)}})^{s'}\oplus 
(P'_{2^{(m)} 1})^{s''},$$ 
where 
$$ s' =  
\left[ \hspace{-0.05in} \begin{array}{c} m-1 \\ k \end{array} \hspace{-0.05in}\right] , 
\quad \quad 
 s'' = 
\left[ \hspace{-0.05in} \begin{array}{c} m-1 \\ k-1 \end{array} \hspace{-0.05in}\right], 
\quad s',s''\in \Z_+[q,q^{-1}] ,  
$$ 
and $P^s$, for a graded module $P$ and $s\in \Z_+[q,q^{-1}]$, denotes 
the direct sum of $s(1)$ copies of $P$ with grading shifts: 
$$P^s = \oplusop{i\in \Z} P^{s_i} \{i \}, \quad \quad s=s(q)=\sum s_i q^i . $$ 
Polynomials $s',s''$ are determined by equation (\ref{eq-five}).  
Applying $\tau$, we obtain 
$$  P_{2^{(k)}12^{(m-k)}} \cong (P_{12^{(m)}})^{s'}\oplus 
(P_{2^{(m)}1})^{s''},$$
so that 
$$[P_{2^{(k)}12^{(m-k)}}] = s' [P_{12^{(m)}}] + s'' [P_{2^{(m)}1}],$$ 
and the equivalent of relation (\ref{eq-five}) holds in $K_0(R)$. 
That (\ref{eq-four}) holds follows from the corresponding result for the nilHecke 
algebra~\cite{KL1, Lau1}.  Therefore, $\gamma$ is well-defined. 
Why $\gamma$ is an isomorphism will be explained next. 
 
\vspace{0.1in} 
     
We have 
\begin{equation} 
 P'_{2^k 1 2^{m-k}}=
R(\alpha_1 + m\alpha_2)1_{2^k 1 2^{m-k}} \ 
\cong \  (P'_{1 2^{(m)}})^{r'} \oplus  (P'_{2^{(m)} 1 })^{{r '}'}
\label{eq-dirsum} 
\end{equation} 
as a graded left $R(\alpha_1 + m \alpha_2)$-module, where 
$$ r' = [k]! [m-k]! s' , \quad \quad r''= [k]! [m]! s'' .$$ 

Projective modules in the above equation correspond to the following three 
idempotents. 

\drawing{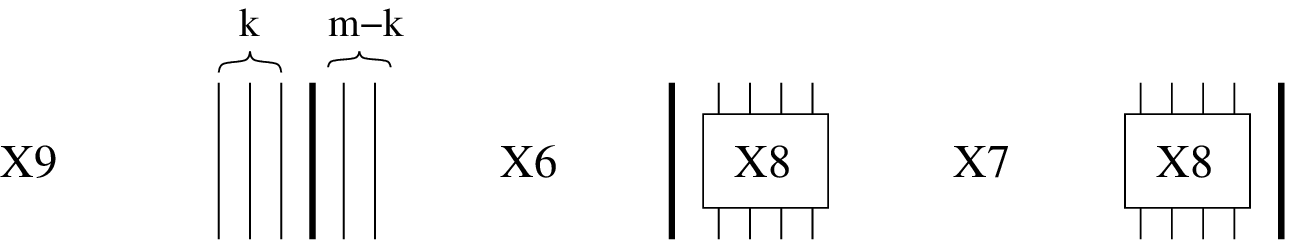} 

Direct sum decomposition (\ref{eq-dirsum}) is equivalent to a choice of 
homogeneous elements 
\begin{eqnarray*} 
\alpha_a' & \in & 1_{2^k 1 2^{m-k}} R(\alpha_1+m\alpha_2) 1_{1 2^{(m)}}, 
\quad \quad 
 1\le a \le r'(1),   \\
\beta_a' & \in & 1_{1 2^{(m)}} R(\alpha_1+m\alpha_2)1_{2^k 1 2^{m-k}}, 
\quad \quad 
1 \le a \le r'(1),   \\
\alpha_a'' & \in &  1_{2^k 1 2^{m-k}} R(\alpha_1+m\alpha_2) 1_{2^{(m)}1} 
\quad \quad 
 1\le a \le r''(1), \\
\beta_a'' & \in & 1_{2^{(m)}1} R(\alpha_1+m\alpha_2)1_{2^k 1 2^{m-k}} 
\quad \quad 
1 \le a \le r''(1),  
\end{eqnarray*} 
such that 
\begin{eqnarray} 
\beta'_b \alpha'_a & = & \delta_{a,b} 1_{1 2^{(m)}}, \\
\beta''_b \alpha''_a & = & \delta_{a,b} 1_{2^{(m)}1}, \\
\beta''_b \alpha'_a & = & 0, \\
\beta'_b \alpha''_a & = & 0, \\    \label{one-decomp} 
1_{2^k 1 2^{m-k}} & = & \sum_{a=1}^{r'(1)} \alpha_a'\beta_a' + 
  \sum_{a=1}^{r''(1)} \alpha_a''\beta_a''. 
\end{eqnarray} 

From the standard representation 
theory of graded rings and the observation that $\mathrm{ker}(\tau)$ belongs 
to the graded Jacobson radical of $R(\alpha_1+m\alpha_2)$ 
 it follows that any finitely-generated graded projective 
$R(1,m)$-module is a direct sum of indecomposable projectives 
$P_{1 2^{(m)}}$ and $P_{2^{(m)}1}$ with grading shifts, that these 
two projectives are not isomorphic, and $K_0(R(1,m))$ is a free 
$\Z[q,q^{-1}]$-module with basis elements $[P_{1 2^{(m)}}]$ and $[P_{2^{(m)}1}]$. 
The above homogeneous elements of $R(\alpha_1+m\alpha_2)$ descend via $\tau$ 
to homogeneous elements of $R(1,m)$  satisfying the same identities. 
We depict these elements by boxes: 

\drawing{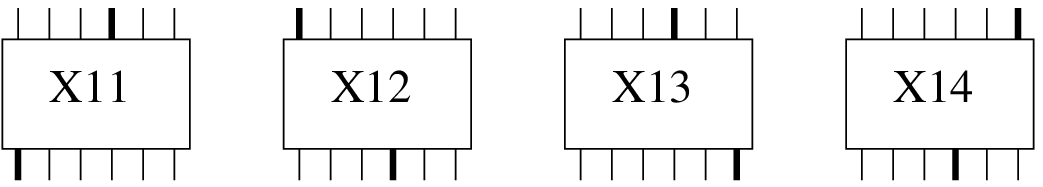} 

Upon applying $\tau$ equation (\ref{one-decomp}) becomes 

\drawing{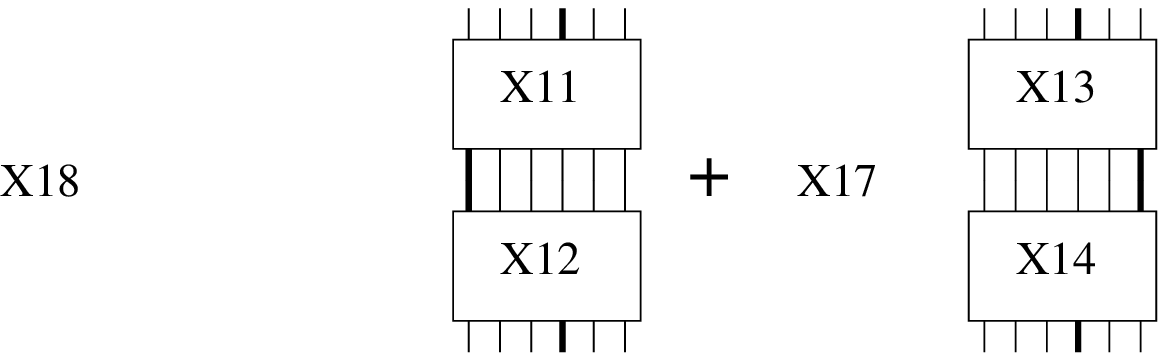} 

To establish that $\gamma$ is an isomorphism, we recall from the beginning of the paper that 
\begin{itemize} 
\item $U^+_{\Z} (0,m)$ is a free $\Z[q,q^{-1}]$-module with generator $E_2^{(m)}$,  
\item $U^+_{\Z} (1,m)$ is a free rank two $\Z[q,q^{-1}]$-module with generators
$E_1 E_2^{(m)}$, $E_2^{(m)} E_1$,  
\item $U^+_{\Z} (2,m)$ is a free $\Z[q,q^{-1}]$-module generated by $E_1E_2^{(m)}E_1$,  
\item $U^+_{\Z} (n,m)=0$ for $n\ge 3$.   
\end{itemize} 

We now look at the size of $K_0(R(n,m))$. 

\vspace{0.06in} 

{\it Case $n=0$.} The ring $R(0,m)$ is graded Noetherian and concentrated in 
cohomological degree $0$ with 
trivial differential. By lemma~\ref{lemma-noeth}, 
$K_0(R(0,m))$ can be computed via finitely-generated 
graded projectives. Any such projective is isomorphic to a finite direct sum of shifts of 
the indecomposable projective $P_{2^{(m)}}$ defined in (\ref{eq-p2m}). Hence,  
$$K_0(R(0,m))  \cong \Z[q,q^{-1}] \cdot [P_{2^{(m)}}],$$ 
giving us a match with $U^+_{\Z} (0,m)$, see above.  

\vspace{0.06in} 

{\it Case $n=1$.} The ring 
$R(1,m)$ is also graded Noetherian concentrated in cohomological degree $0$ with 
trivial differential. We have established that the graded 
Grothendieck group $K_0(R(1,m))$ is a free rank two $\Z[q,q^{-1}]$-module 
with the basis given by symbols $[P_{12^{(m)}}]$, $[P_{2^{(m)}1}]$. This matches
with the above basis for $U^+_{\Z} (1,m)$. 

\vspace{0.06in} 

{\it Case $n\ge 3$.} We will prove that dg algebra $R(n,m)$ has trivial homology 
when $n\ge 3$. Consider the element $y_k \in R(3,m)$: 

\drawing{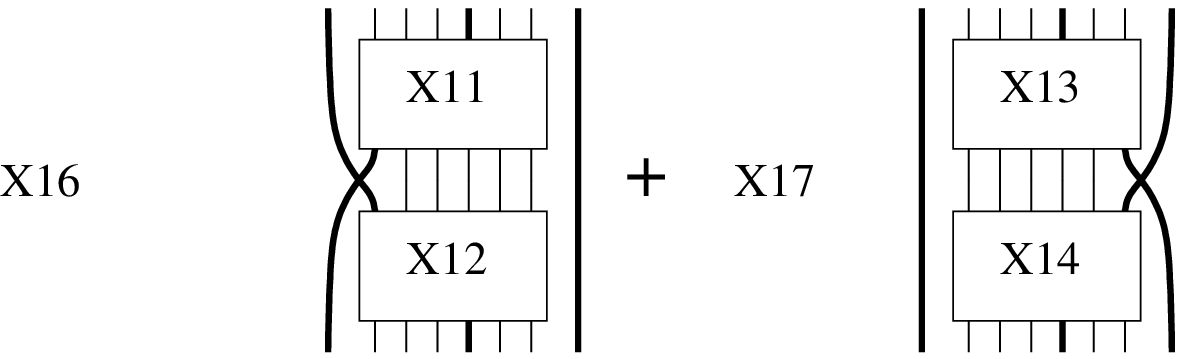}  

Applying $d$ and using the previous diagrammatic equation, we get
$$d(y_k) = 1_{1 2^k 1 2^{m-k} 1}.$$ 
Therefore, for any sequence ${\bf i}$ which 
contains at least three ones we can find $y_{\bf i}\in R(n,m)$ such that $d(y_{\bf i})=1_{\bf i}$. 
We can write $1\in R(n,m)$ as the sum of $1_{\bf i}$ over all sequences 
${\bf i}$ with $n$ ones ($n\ge 3$) and $m$ twos. Then 
$$ d(\sum_{\bf i} y_{\bf i}) =1 $$ 
implying that $H(R(n,m))=0$ and $K_0(R(n,m))=0$. 

\vspace{0.06in} 

{\it Case $n=2$.} Observe that the unit element of $R(2,m)$ decomposes 
$$ 1 = \sum_{k,\ell} 1_{2^k 1 2^{\ell} 1 2^{m-k -\ell}}$$
and that projective module $P_{2^k 1 2^{\ell} 1 2^{m-k -\ell}}$ is isomorphic 
to the direct sum of shifts of modules $P_{1 2^{(k + \ell)} 1 2^{m-k -\ell}}$ and
$P_{2^{(k + \ell)} 1 1 2^{m-k -\ell}}$. The latter are 
contractible, as complexes of vector spaces (since their sequences contain two consecutive 
ones), while modules of the first kind 
decompose as direct sum of shifts of modules $P_{11 2^{(m)}}$ and $P_{1 2^{(m)}1}$. 
We can discard $ P_{11 2^{(m)}}$ because of contractibility and conclude that 
the $K_0$ groups of graded dg rings $R(2,m)$ and 
$$ R_1 \ := \ 1_{1 2^{(m)} 1} R(2,m) 1_{1 2^{(m)} 1} $$ 
are naturally isomorphic. To make this argument more accurate, introduce 
additional dg rings defined via idempotents $e(4),e(3), e(2)$: 
\begin{eqnarray*} 
R_4 & = & e(4) R(2,m) e(4) , \quad e(4) \ :=\  \sum_{t=0}^m (1_{1 2^{(t)} 1 2^{m-t}}+ 
  1_{2^{(t)} 11 2^{m-t}}), \\
R_3 & = & e(3) R_4 e(3), \quad e(3) \ := \ \sum_{t=0}^m 1_{1 2^{(t)} 1 2^{m-t}}, \\
R_2 & = &  e(2) R_3 e(2), \quad e(2) \ := \ 1_{11 2^{(m)}} + 1_{1 2^{(m)} 1}. 
\end{eqnarray*} 
Graded dg rings $R_4$ and $R(2,m)$ are graded dg Morita equivalent in the strongest 
sense, via dg bimodules $R(2,m)e(4)$ and $e(4)R(2,m)$. These bimodules produce an 
equivalence of abelian categories of graded dg modules over $R_4$ and $R(2,m)$, as well as 
all the other categories of dg modules and their graded counterparts that appear 
in the diagram (\ref{diag-cd}). 
Hence, there is a canonical isomorphism $K_0(R(2,m))\cong K_0(R_4)$. 
Next, $R_3$ is a graded dg subring of $R_4$ obtained by removing 
chunks of $R_4$ corresponding to contractible  idempotents $1_{2^{(t)} 11 2^{(m-t)}}$. 
In view of lemma~\ref{throw-idemp} and invariance of $K_0$ under 
quasi-isomorphisms of graded dg rings there is a natural isomorphism 
$K_0(R_4)\cong K_0(R_3)$.  Graded dg rings $R_3$ and $R_2$ are graded 
dg Morita equivalent, via bimodules 
$$ R_3 (1_{11 2^{(m)}} + 1_{1 2^{(m)} 1}) \quad \textrm{and} \quad 
(1_{11 2^{(m)}} + 1_{1 2^{(m)} 1})R_3, $$
giving us an isomorphism $K_0(R_3) \cong K_0(R_2)$. Finally, 
reducing $R_2$ via lemma~\ref{throw-idemp} applied to the contractible idempotent 
$1_{11 2^{(m)}}$ results in $R_1$.  Nonunital inclusions 
$$ R_1 \subset R_2 \subset R_3 \subset R_4 \subset R(2,m)$$ 
are quasi-isomorphisms, giving rise to canonical isomorphisms of 
Grothendieck groups 
\begin{eqnarray} \label{eq-congs1}
 K_0(R_1) \cong K_0(R_2) \cong K_0(R_3) \cong K_0(R_4) \cong K_0(R(2,m)). 
\end{eqnarray} 
Diagrams representing elements of $R_1$ 
have two fermionic lines that start and end at the leftmost and rightmost top and bottom 
endpoints and $m$ bosonic lines, capped off on both sides by the projector $e_m=1_{2^{(m)}}$.  
We can decompose $R_1$ into the direct sum of the 2-sided ideal $I$ spanned by 
diagrams in which the two fermionic lines intersect and the subring $R_1'$ spanned by 
diagrams with the disjoint fermionic lines, $R_1 = R_1'\oplus I$. The ring 
$R_1' $ is isomorphic to $e_m H_m e_m $, where $H_m$ is the nilHecke algebra
(note that  nonintersecting fermionic lines can be pulled away to be disjoint from 
the $m$ bosonic lines). In turn, $e_m H_m e_m $ is isomorphic to the 
ring of symmetric polynomials in $x_1, \dots, x_m$, corresponding to the dots 
on the bosonic lines, so that $R_1'\cong Sym(x_1, \dots, x_m)$. 
Ideal $I$ is isomorphic to $R_1'$ when viewed as an $R_1'$-bimodule, with 
the generator $X$ depicted below and on the left. The 
differential in the dg ring $R_1$ is zero on $R_1'$ and takes $I$ injectively to $R_1'$, 
$$ 0 \lra I \stackrel{d}{\lra} R_1' \lra 0, $$ 
with the generator $X$  taken to $x_1 x_2 \dots x_m$. 
In the diagram below each box denotes idempotent $e_m$. 

\drawing{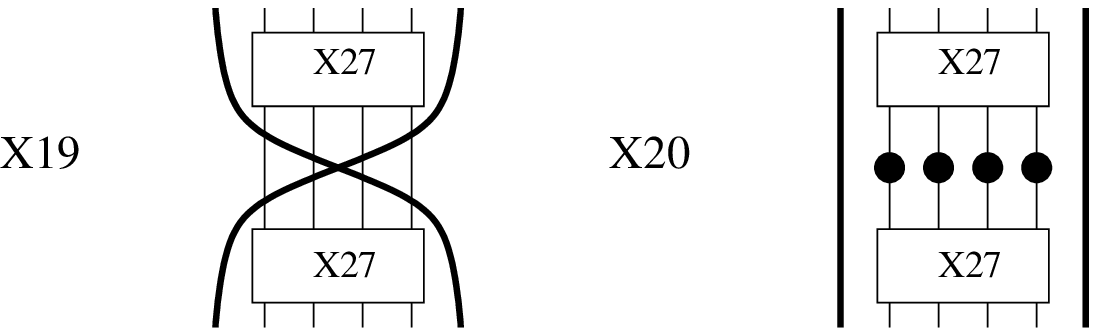} 

Therefore, the quotient map $R_1 \lra R_1'/d(I)$ is a quasi-isomorphism of 
dg rings, inducing an isomorphism $K_0(R_1) \cong K_0(R_1'/d(I))$.  
The ring  $R_1'/d(I)$  
has trivial differential and is concentrated in cohomological degree $0$. 
It is isomorphic to the quotient of the graded ring of symmetric 
polynomials $Sym(x_1, \dots, x_m)$ by the ideal $(x_1x_2\dots x_m)$. 
This quotient is a graded local ring, with $K_0$ isomorphic to $\Z[q,q^{-1}]$. 
The above arguments tell us that 
$$ K_0(R(2,m)) \ \cong \ K_0(R_1) \ \cong \  \Z[q,q^{-1}],$$ 
with the generator being $[P_{1 2^{(m)} 1}]$. 

\vspace{0.05in} 

We conclude that the map $\gamma$ is an isomorphism. It is straightforward 
to check that $\gamma$ respects comultiplication of twisted bialgebras 
in (\ref{eq-gamma}). Hence, we have a canonical isomorphism of twisted bialgebras 
$$  U^+_{\Z}  \ \cong \ K_0(R) $$ 
taking weight spaces of $U^+_{\Z}$ to $K_0(R(n,m))$ for 
various $n,m$. 

\vspace{0.15in} 
 

\noindent
{\bf Perspectives} 


\noindent 
It is not hard to guess how one should couple the above diagrammatics
 for categorified $U^+_q(\mf{gl}(1|2))$ to the diagrammatics~\cite{KL1}  
for categorified $U^+_q(\mf{g})$ to produce a graphical calculus for 
categorified $U^+_q(\mf{gl}(1|n))$ and some other classical Lie superalgebras 
in place of $\mf{gl}(1|n)$. 
Technical obstacles related to switching from algebras to dg algebras and 
their representations were largely avoided in the present paper due to small size of 
Grothendieck groups of dg rings $R(n,m)$ and manual case-by-case considerations. 
We plan to discuss categorification of $U^+_q(\mf{gl}(1|n))$ in a follow-up paper. 
Categorification of $U^+_q(\mf{gl}(n|m))$ for $n,m\ge 2$ should require additional 
ideas. 

Varagnolo and Vasserot~\cite{VV} established that rings $R(\nu)$ from~\cite{KL1},  
\cite{Rou}  
that categorify weight spaces of $U^+_q(\mf{g})$ are isomorphic to equivariant ext  
groups of sheaves on Lusztig quiver varieties. It is an interesting problem 
to find a similar interpretation for the dg rings $R(n,m)$ described above that 
categorify weight spaces of $U^+_q(\mf{gl}(1|2))$ and for their generalizations.   
  
The problem of categorifying the entire (suitably idempotented) quantum
group $U_q(\mf{gl}(1|1))$ rather than just its positive half is open as well. With Lauda's
categorification~\cite{Lau1} of $U_q(\mf{sl}(2))$ in mind, we expect that
one should generalize the Lipshitz-Ozsv\'ath-Thurston rings  to allow the
lines to travel in all directions rather than just up. 
Another challenging problem is to extend the work of 
Webster~\cite{W1, W2} to the superalgebra case. 

\vspace{0.3in} 


\noindent 
{\bf Acknowledgments} 


\noindent 
The author would like to thank the NSF for partial support via 
grant DMS-0706924.

\vspace{0.3in} 

\noindent 
{\bf References} 


\vspace{-0.35in} 

\def\refname{}

\vspace{0.2in} 

{ \sl \small Department of Mathematics, Columbia University, New York, NY 10027}
 
{  \tt \small email: khovanov@math.columbia.edu}


\begin{thebibliography}{mmmm}

\bibitem[1]{BL} J.~Bernstein and  V.~Lunts, 
{\sl Equivariant sheaves and functors}, LNM 1578, (1994). 

\bibitem[2]{Kac1} V.~G.~Kac, {\sl Lie superalgebras}, Adv. in Math. {\bf 26} (1977), 
8--96. 

\bibitem[3]{Kh1} 
M.Khovanov, {\sl NilCoxeter algebras categorify the Weyl algebra}, Comm. Algebra, {\bf 29} (11) 
5033--5052, 2001, math.RT/9906166. 

\bibitem[4]{KL1}
M.~Khovanov and A.~Lauda, {\sl A diagrammatic approach to categorification of quantum groups I},
Represent. Theory {\bf 13} (2009), 309--347,  arxiv 0803.4121.

\bibitem[5]{Lau1}
A.~Lauda, {\sl A categorification of quantum sl(2)}, to appear in Advances in Math., 
arxiv 0803.3848. 

\bibitem[6]{LOT1} 
R.~Lipshitz, P.~Ozsv\'ath, and D.~Thurston, {\sl Bordered Heegaard Floer 
homology: invariance and pairing}, arxiv 0810.0687. 

\bibitem[7]{LOT2} 
R.~Lipshitz, P.~Ozsv\'ath, and D.~Thurston, {\sl Slicing planar grid diagrams: a gentle 
introduction to bordered Heegaard Floer homology}, arxiv 0810.0695. 
 
\bibitem[8]{Lus1} 
G.~Lusztig, {\sl Introduction to quantum groups}, vol.~110 of \emph{Progress in Mathematics}, 
Birkh\"auser Boston, 1993. 

\bibitem[9]{MOS} 
C.~Manolescu, P.~Ozsv\'ath, and S.~Sarkar, {\sl A combinatorial description 
of knot Floer homology}, Annals of Math. {\bf 169} (2009), no.2, 633--660, math/0607691. 

\bibitem[10]{Rou} 
R.~Rouquier, {\sl 2-Kac-Moody algebras}, arxiv 0812.5023. 
 
\bibitem[11]{VV} 
M.~Varagnolo and E.~Vasserot, {\sl Canonical bases and Khovanov-Lauda algebras}, arxiv 0901.3992. 

\bibitem[12]{W1}
B.~Webster, {\sl Knot invariants and higher representation theory I: diagrammatic and geometric categorification of 
tensor products}, arXiv 1001.2020. 

\bibitem[13]{W2}
B.~Webster, {\sl Knot invariants and higher representation theory II: the categorification of quantum knot invariants}, 
arXiv 1005.4559. 

\end{thebibliography}
\end{document}